\documentclass{article}
\usepackage{amsmath}
\usepackage{amssymb, amscd, theorem, amsfonts, latexsym}
\usepackage{enumerate}
\usepackage[dvips]{graphicx,color}

\newcommand{\no}{\noindent}
\newtheorem{thm}{Theorem}[section]
\newtheorem{prop}[thm]{Proposition}
\newtheorem{cor}[thm]{Corollary}
\newtheorem{rem}[thm]{Remark}
\newtheorem{definition}[thm]{Definition}
\newtheorem{ex}[thm]{Example}
\newtheorem{lem}[thm]{Lemma}

\newcommand{\R}{\mathbb{R}}

\newcommand{\pd}{\partial}
\newcommand{\al}{\alpha}
\renewcommand{\phi}{\varphi}
\renewcommand{\epsilon}{\varepsilon}

\begin{document}

\title{\bf Extremal problems for the solid angle}

\author{\bf Shigehiro Sakata}

\date{\today}

\maketitle

\def\vect#1{\mbox{\boldmath $#1$}}
\def\vecttiny#1{\mbox{\tiny \boldmath $#1$}} 

\begin{abstract}
We consider the solid angle that a planar compact subset subtends at a point in a level set of height $h$ and study two extremal problems for the solid angle. One of the variables is a point in such a plane, that is, we study the properties of the solid angle maximizer. The other is the pair of a planar compact subset and a point.  \\

\hspace{-1.5em}{\it keywords and phrases}. solid angle, barycenter, extremal problem, $r^{\al - m}$-potential, distance kernel, moving plane method.
\end{abstract}
%%%%%%%%%%%%%%%%%%%%%%%%%%%%%%%%%%%%%%%%%%%%%%%%%%%%%%%%%%%%%%%%%%%%%%%%%%%%%%%%%%%%%%%%%%%%%%%%%%%%%%%%%%%%%%
\section{Introduction}

In [5] Katsuyuki Shibata introduced a {\it spatial illuminating center} of a triangle $\triangle$ of height $h>0$, which is a point that maximizes the {\it ``total brightness''} of a triangular park $\triangle$ if a light source is located above that point by height $h$. It is easy to see that the total brightness is proportional to the solid angle that $\triangle$ subtends at the point of the light source if the height $h$ is fixed. Let us generalize the notion of spatial illuminating center for compact subsets $\Omega$ in $\R^2$ with non-empty interior, and call it the {\it solid angle center} in what follows.\\
\hspace{1.5em}The existence of the solid angle center is rather straightforward, but the uniqueness does not always hold. In fact, a disjoint union of two discs of the same size has (at least) two solid angle centers of height $h>0$ if $h$ is enough small. We study the condition for the uniqueness of the solid angle center.\\
\hspace{1.5em}This type of research originated in [3] by Jun O'Hara, where he introduced renormalization of the $r^{\al -m}$-potential of the characteristic function of a compact domain in $\R^m$ for $\al \leq 0$ (we remark that the $r^{\al -m}$-potential is called the Riesz potential when $\al$ satisfies $0< \al <m$ and is called the logarithmic potential when $\al=m=2$.), and showed the uniqueness of a point that attains the maximum value of the (renormalizaiton of) $r^{\al -m}$-potential if $\al \leq 1$ and the domain is convex. On the other hand, it is easy to show that for any compact domain there is a unique point that minimizes the $r^{\al -m}$-potential if $\al \geq m+1$. Thus he obtained a generalization of the barycenter, which corresponds to the case when $\al =m+2$.\\
\hspace{1.5em}The main idea of [3] is to show that if the condition above is satisfied then the second derivative of (the renormalization of) the potential is negative on the domain using the boundary integral expression. Unfortunately, this does not hold as it did for our case. To be precise, the second derivative of the solid angle function is not always negative on the domain even if the domain is convex.\\
\hspace{1.5em}In this paper we study some properties of the solid angle center and the solid angle function, in particular the uniqueness of the solid angle center for the case when the above method does not work.\\
\hspace{1.5em}We show the following properties:
\begin{enumerate}
\item[(a)] The solid angle center of $\Omega$ of height $h$ converges to the barycenter of $\Omega$ as $h$ goes to $+\infty$. 
\item[(b)] The maximum value of the solid angle function of $\Omega$ is not greater than that of a disc with the same area. The equality holds if and only if $\Omega$ is a disc.
\item[(c)] The integral of the solid angle function over $\Omega$ is not greater than that of a disc with the same area. The equality holds if and only if $\Omega$ is a disc.
\item[(d)] If $\Omega$ is convex, then the solid angle fucntion is superharmonic on $\Omega$.
\item[(e)] The solid angle center of $\Omega$ of height $h>0$ is unique if one of the following conditions is satisfied. 
\begin{enumerate}
\item[(1)] If $h$ is not smaller than twice the diameter of $\Omega$. An improved condition will be given later. 
\item[(2)] If $\Omega$ is convex and axially symmetric. 
\item[(3)] If $\Omega$ is convex and $h$ is smaller than a constant that will be given later.
\end{enumerate}
\end{enumerate}
We remark that $\Omega$ is an arbitrary compact planar subset with non-empty interior in the case (1) above and the height $h$ is arbitrary in the case (2). It turns out that the same argument for the second condition also works for the $r^{\al -2}$-potential (section 5).\\
\hspace{1.5em} Our results for the uniqueness of the solid angle center are obtained by the following procedure. Using radial symmetry of the potential and moving plane method (cf. [1]), we can show that there is a smaller region, depending on $\Omega$, so that the complement has no chance to have any center on in. This region was introduced in [3] and is called the {\it minimal unfolded region} of $\Omega$.
\begin{description}
\item[Step1.] Compute the second derivative of the solid angle function for any direction.
\begin{description}
\item[(1-1)] For the first condition in (e) above, the computation is direct.
\item[(1-2)] For the second condition, we use the contour integral expression.
\item[(1-3)] For the third condition, we use the radial function.
\end{description}
\item[Step2.] Estimate the signature of the second derivative on the minimal unfolded region of $\Omega$.
\item[Step3.] If $\Omega$ has two solid angle centers, then we obtain a contradiction to the statement of step 2.
\end{description}

{\bf Acknowledment.}
The author would like to express his deep gratitude to his advisor Jun O'Hara for informing him of these problems, the uniqueness of the solid angle center and the $r^{\al -m}$-center, and for giving many helpful advices to him.
%%%%%%%%%%%%%%%%%%%%%%%%%%%%%%%%%%%%%%%%%%%%%%%%%%%%%%%%%%%%%%%%%%%%%%%%%%%%%%%%%%%%%%%%%%%%%%%%%%%%%%%%%%%%%%
\section{Preliminaries}
%%%%%%%%%%%%%%%%%%%%%%%%%%%%%%%%%%%%%%%%%%%%%%%%%%%%%%%%%%%%%%%%%%%%%%%%%%%%%%%%%%%%%%%%%%%%%%%%%%%%%%%%%%%%%%
\subsection{Notations and definition}
Let $\Omega$ be a compact subset in $\R^2$ with non-empty interior. Let $x$ be a point in $\R^2$ and $h$ be a fixed positive number. Define the map $p_x^{(h)}:\Omega \to S^2$ by
\[
p_x^{(h)}(y) = \frac{(y,0)-(x,h)}{\left\vert (y,0)-(x,h) \right\vert} = \frac{(y-x,0-h)}{\sqrt{r^2+h^2}},
\] where $r= \lvert y - x \rvert$. The solid angle at a point $(x,h)$, denoted by $A_\Omega^{(h)}(x)$, is given by the area of the image of $p_x^{(h)}$. By direct calculation, we see that $A_\Omega^{(h)}(x)$ is given by 
\[
A_\Omega^{(h)}(x) = \int_\Omega \frac{h}{(r^2+h^2)^{3/2}}dy.
\] 
\begin{definition}\label{def_solid}{\rm A point $x$ in $\R^2$ is the {\it solid angle center} of $\Omega$ of height $h>0$ if it attains the maximum value of $A_\Omega^{(h)}: \R^2 \to \R$.}
\end{definition}

We denote a point in $\R^2$ by $x=(x_1, x_2)$ and a point in $\Omega$ by $y=(y_1,y_2)$. We understand that a letter $r$ is always used for $r=\lvert y  - x \rvert$.
%%%%%%%%%%%%%%%%%%%%%%%%%%%%%%%%%%%%%%%%%%%%%%%%%%%%%%%%%%%%%%%%%%%%%%%%%%%%%%%%%%%%%%%%%%%%%%%%%%%%%%%%%%%%%%
\subsection{Minimal unfolded region and existence of the solid angle center}
We show the existence of the solid angle center using the same argument as in [3]. We denote the convex hull of $\Omega$ by ${\rm conv}(\Omega)$.

\begin{prop}\label{prop_exist} There exists a solid angle center of $\Omega$ in ${\rm conv}(\Omega)$.
\end{prop}
{\bf Proof.} Since $A_\Omega^{(h)}$ is continuous on $\R^2$ and ${\rm conv}(\Omega)$ is compact, there exists a point $x_c$ that attains the maximum value of $A_\Omega^{(h)}$ in ${\rm conv}(\Omega)$. We show that $x_c$ is a solid angle center of $\Omega$.\\
\hspace{1.5em}Let $x$ be a point in $({\rm conv}(\Omega))^c$. Let $x'$ be a point on $\pd {\rm conv}(\Omega)$ such that ${\rm dist} (x,\pd {\rm conv}(\Omega)) = \lvert x-x' \rvert$. Then ${\rm conv}(\Omega)$ is contained in a half-space whose boundary is the line orthogonal to a line through $x$ and $x'$. Therefore, for any point $y$ in $\Omega$, we have $\lvert x-y\rvert > \lvert x'-y \rvert$.\\
\hspace{1.5em}Hence we obtain $A_\Omega^{(h)}(x)<A_\Omega^{(h)}(x') \leq A_\Omega^{(h)}(x_c)$. \hspace{\fill}$\Box$\\

In order to improve Proposition \ref{prop_exist}, we introduce some results in [3] with a slight modification.

\begin{definition}\label{def_folded}{\rm ([3]) Let $v$ be a point in the unit sphere $S^1$. Let $I_{v,c} \ \ (c \in \R)$ be a reflection of $\R^2$ in a line $\left\{ z \in \R^2 \vert z \cdot v =c \right\}$. Put 
\[
\Omega_{v,a}^+ = \Omega \cap \left\{ z \in \R^2 \vert z \cdot v \geq a \right\} \ \
\] and
\[
\Omega_{v,a}^- = \Omega \cap \left\{ z \in \R^2 \vert z \cdot v \leq a \right\}.
\] Let $u(v)$ and $l(v)$ be given by 
\begin{align*}
u(v) &= \inf \left\{ b \in \R \vert I_{v,c}(\Omega_{v,c}^+) \subset \Omega, \ \ \forall c \geq b \right\},\\
l(v) &= \sup \left\{ b \in \R \vert I_{v,c}(\Omega_{v,c}^-) \subset \Omega, \ \ \forall c \leq b \right\}.
\end{align*}
Define the {\it minimal unfolded region} of $\Omega$ by
\[
uf(\Omega) = \bigcap_{v \in S^1} \left\{ z \in \R^2 \vert l(v) \leq z \cdot v \leq u(v) \right\}.
\]}
\end{definition}

\begin{rem}\label{rem_folded}{\rm ([3]) $uf(\Omega)$ is compact and contained in ${\rm conv}(\Omega)$.}
\end{rem}

\begin{thm}\label{thm_folded} For any $h>0$ a solid angle center of $\Omega$ belongs to $uf(\Omega)$.
\end{thm}
{\bf Proof.} We use the radial symmetry of partial derivatives of $A_\Omega^{(h)}$ and the moving plane method ([1]) in parallel to that of Theorem 4.6 in [3].\\
\hspace{1.5em} Suppose a point $x$ is not contained in $uf(\Omega)$. As partial derivative of $A_\Omega^{(h)}$ have symmetry,
\[
\frac{\pd A_{\Omega}^{(h)}}{\pd x_1}(x) = \frac{\pd A_{\Omega_1}^{(h)}}{\pd x_1}(x) +\frac{\pd A_{\Omega_2}^{(h)}}{\pd x_1}(x) =\frac{\pd A_{\Omega_2}^{(h)}}{\pd x_1}(x) = 3 \int_{\Omega_2} \frac{y_1 - x_1}{(r^2+h^2)^{5/2}} dy \neq 0,
\] which is contradiction (see figure 1). \hspace{\fill}$\Box$\\
%%%%%%%%%%%%%%%%%%%%%%%%%%%%%%%%%%%%%%%%%%%%%%%%%%%%%%%%%%%%%%%%%%%%%%%%%%%%%%%%%%%%%%%%%%%%%%%%%%%%%%%%%%%%
\begin{figure}[hbtp]
\centering
\scalebox{0.40}{\includegraphics[clip]{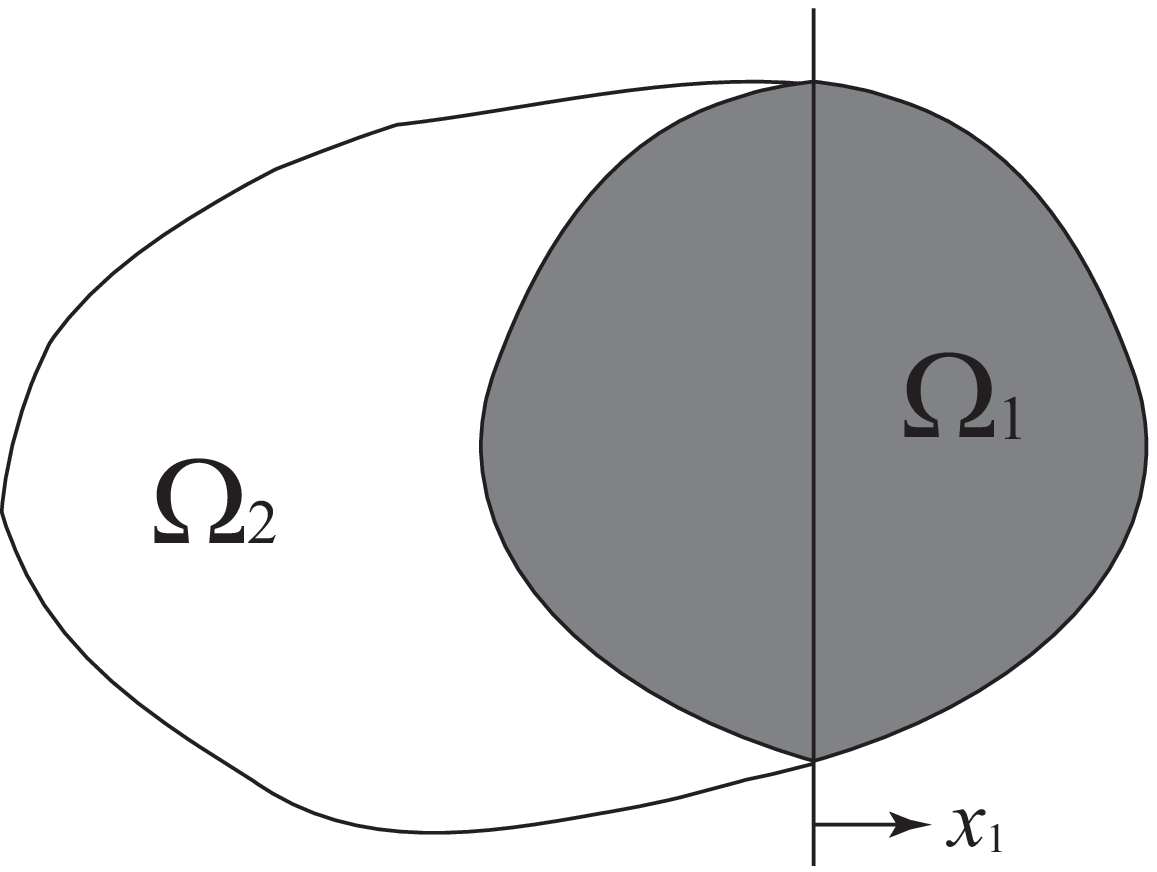}}
\caption{$\displaystyle \frac{\pd A_{\Omega}^{(h)}}{\pd x_1}(x) = 3 \int_{\Omega_2} \frac{y_1-x_1}{(r^2+h^2)^{5/2}} dy \neq 0$}
\label{moving_plane}
\end{figure}
%%%%%%%%%%%%%%%%%%%%%%%%%%%%%%%%%%%%%%%%%%%%%%%%%%%%%%%%%%%%%%%%%%%%%%%%%%%%%%%%%%%%%%%%%%%%%%%%%%%%%%%%%%%%
\section{Some properties}
%%%%%%%%%%%%%%%%%%%%%%%%%%%%%%%%%%%%%%%%%%%%%%%%%%%%%%%%%%%%%%%%%%%%%%%%%%%%%%%%%%%%%%%%%%%%%%%%%%%%%%%%%%%%
\subsection{Limit point of the solid  angle center as {\boldmath $h$} goes to {\boldmath $+\infty$}}
\begin{rem}\label{barycenter} {\rm ([3]) The barycenter of $\Omega$ is given by the following coordinate and we can obtain the barycenter as the minimum point of a map $\displaystyle \R^2 \ni x \mapsto \int_\Omega r^2 dy \in \R$:
\[
x = \frac{1}{{\rm Area}(\Omega)} \int_\Omega y dy.
\]} 
\end{rem}
\begin{thm}\label{limit_point} The solid angle center of $\Omega$ of height $h$ converges to the barycenter of $\Omega$ as $h$ goes to $+\infty$.
\end{thm}
{\bf Proof.} Let $G$ be the barycenter of $\Omega$ and $D_\epsilon(G)$ the $\epsilon$-neighborhood of $G$. Note that
\[
A_\Omega^{(h)} (x) = \frac{1}{h^2} \int_\Omega \left( -\frac{3}{2} \left( \frac{r}{h} \right)^2 + \left( 1+ \left( \frac{r}{h} \right)^2 \right)^{-\frac{3}{2}} + \frac{3}{2} \left( \frac{r}{h} \right)^2 \right) dy.
\] For any $x \in {\rm conv}(\Omega)$, we have
\begin{align*}
\left\vert \frac{\pd}{\pd x_j}  \int_\Omega \left(\left( 1+ \left( \frac{r}{h} \right)^2 \right)^{-\frac{3}{2}} + \frac{3}{2} \left( \frac{r}{h} \right)^2 \right) dy \right\vert
&\leq \frac{1}{h} \left\vert \int_\Omega \left( -3 \left( 1+ \left( \frac{r}{h} \right)^2 \right)^{-\frac{5}{2}} \frac{r}{h} + 3\frac{r}{h} \right) dy \right\vert \\
&\leq \frac{3 {\rm diam}(\Omega)}{h^2} \int_\Omega \left( 1 - \left( 1+ \left( \frac{r}{h} \right)^2 \right)^{-\frac{5}{2}} \right) dy \\
&\leq \frac{3 {\rm diam}(\Omega)}{h^2} \int_\Omega \frac{5}{2} \left( \frac{r}{h} \right) ^2 dy \\
&\leq \frac{15 \left( {\rm diam} (\Omega ) \right)^3 {\rm Area}(\Omega)}{2h^4}.
\end{align*} Put $\displaystyle M = \frac{15}{2} \left( {\rm diam} (\Omega ) \right)^3 {\rm Area}(\Omega)$. Fix an arbitrary $\epsilon >0$, and let
\[
m(\epsilon) = \max_{j \in \{1,2\}} \left(\inf_{x \in {\rm conv}(\Omega) \backslash D_\epsilon (G)} \frac{3}{2} \left\vert \frac{\pd}{\pd x_j} \int_\Omega r^2 dy \right\vert \right) >0.
\] If $\displaystyle h \geq \sqrt{\frac{M}{m(\epsilon)}}$, then for any $x \in {\rm conv}(\Omega) \backslash D_\epsilon (G)$, we have
\[
\max_{j \in \{1, 2\}} \left( \inf_{x \in {\rm conv}(\Omega) \backslash D_\epsilon (G)} \left\vert \frac{\pd A_\Omega^{(h)}}{\pd x_j} (x) \right\vert \right) \geq \frac{1}{h^4} \left( m(\epsilon ) - \frac{M}{h^2} \right) >0.
\]
\hspace{1.5em}Hence we obtain the conclusion.\hspace{\fill}$\Box$
%%%%%%%%%%%%%%%%%%%%%%%%%%%%%%%%%%%%%%%%%%%%%%%%%%%%%%%%%%%%%%%%%%%%%%%%%%%%%%%%%%%%%%%%%%%%%%%%%%%%%%%%%%%%
\subsection{Extremal problems for {\boldmath $A_\Omega^{(h)}$}}
We study two extremal problems. One of the functionals is the maximum value of $A_\Omega^{(h)}$ and the other is the integral of $A_\Omega^{(h)}$ over $\Omega$. Let $D$ be a disc in $\R^2$. By the same argument as in [4] we have the following.
\begin{thm}\label{variation} If ${\rm Area}(\Omega) = {\rm Area}(D)$, then 
\[
\max_{x \in \R^2} A_\Omega^{(h)} (x) \leq \max_{x \in \R^2} A_D^{(h)} (x)
\] and that equality holds if and only if $\Omega$ is a disc.
\end{thm}
{\bf Proof.} Theorem \ref{thm_folded} guarantees that the center of $D$ is the solid angle center of D. By a translation of $\R^2$, we may assume that the center of $D$ and the solid angle center of $\Omega$ coincide with the origin $0$. Since $\left\vert y' \right\vert \leq \left\vert y'' \right\vert$ for $y' \in D \backslash (D \cap \Omega)$ and $y'' \in \Omega \backslash (D \cap \Omega)$, we obtain
\[
A_D^{(h)} (0) - A_\Omega^{(h)} (0) = \int_{D \backslash (D \cap \Omega)} \frac{h}{\left( \left\vert y \right\vert ^2 +h^2 \right)^{3/2}} dy -\int_{\Omega \backslash (D \cap \Omega)} \frac{h}{\left( \left\vert y \right\vert ^2 +h^2 \right)^{3/2}} dy \geq 0.
\] That equality holds if and only if ${\rm Area}(\Omega \backslash (D \cap \Omega))=0$, that is, $\Omega$ is a disc. \hspace{\fill}$\Box$
\begin{thm}\label{Morgan} $([2])$ Let $f(t)$ be a strictly decreasing positive function defined for all $t>0$ with $f(r)$ locally integrable over $\R^2 \times \R^2$. If ${\rm Area}(\Omega)={\rm Area}(D)$, then 
\[
\int_{\Omega \times \Omega} f(r) dxdy \leq \int_{D \times D} f(r) dxdy
\] and that equality holds if and only if $\Omega$ is a disc.
\end{thm}
\begin{cor}\label{variation_integral} If ${\rm Area}(\Omega) = {\rm Area}(D)$, then 
\[
\int_\Omega A_\Omega^{(h)} (x)dx \leq \int_D A_D^{(h)} (x)dx
\] and that equality holds if and only if $\Omega$ is a disc.
\end{cor}
{\bf Proof.} Put $\displaystyle f(r) = \frac{1}{\left( r^2 +h^2 \right)^{3/2}}$. \hspace{\fill}$\Box$\\
%%%%%%%%%%%%%%%%%%%%%%%%%%%%%%%%%%%%%%%%%%%%%%%%%%%%%%%%%%%%%%%%%%%%%%%%%%%%%%%%%%%%%%%%%%%%%%%%%%%%%%%%%%%%%
\subsection{Laplacian}
\begin{prop}\label{prop_contour} Suppose that $\Omega$ has a piecewise $C^1$ boundary. Then $\displaystyle \frac{\pd A_\Omega^{(h)}}{\pd x_j}$ can be expressed by the contour integral on $\pd \Omega$ as
\begin{equation}\label{contour}
\frac{\pd A_\Omega^{(h)}}{\pd x_j}(x) = -h \int_{\pd \Omega} \frac{e_j \cdot n(y)}{(r^2+h^2)^{3/2}} ds,
\end{equation} where $e_j$ is the $j$-th unit vector of $\R^2$, $n(y)$ is the unit outer noromal to $\pd \Omega$ at $y$ and $s$ denotes the arc-length parameter of $\pd \Omega$.
\end{prop}
{\bf Proof.} The proof is parallel to that of Proposition 3.2 in [3]. Note that
\[
\frac{\pd}{\pd x_j}\left( \frac{1}{(r^2+h^2)^{3/2}}\right) =-\frac{\pd}{\pd y_j}\left( \frac{1}{(r^2+h^2)^{3/2}}\right) 
\] and
\[
 \int_{\Omega}\frac{\pd}{\pd y_j} \left(\frac{1}{(r^2+h^2)^{3/2}} \right) dy = \left\{
\begin{array}{ll}
\displaystyle    \int_{\pd \Omega}\frac{1}{(r^2+h^2)^{3/2}} dy_2  &({\rm if} \ \ j=1), \\
\displaystyle -  \int_{\pd \Omega}\frac{1}{(r^2+h^2)^{3/2}} dy_1  &({\rm if} \ \ j=2), \\
\end{array}
\right. 
\] by Stokes' theorem. We can see the conclusion by using these remarks. \hspace{\fill}$\Box$
\begin{thm}\label{laplacian} If $\Omega$ is convex and has a piecewise $C^1$ boundary. Then $A_{\Omega}^{(h)}$ is superharmonic on $\Omega$.
\end{thm}
{\bf Proof.} Note that, for any $x \in \Omega$ and $y \in \pd \Omega$, $(y - x) \cdot n(y) \geq 0$ since $\Omega$ is convex. It follows from Proposition \ref{prop_contour} that 
\[
\triangle A_{\Omega}^{(h)} (x) = -3h \int_{\pd \Omega} \frac{(y - x)\cdot n(y)}{(r^2 + h^2)^{5/2}}ds <0,
\] for $x \in \Omega$. \hspace{\fill}$\Box$
%%%%%%%%%%%%%%%%%%%%%%%%%%%%%%%%%%%%%%%%%%%%%%%%%%%%%%%%%%%%%%%%%%%%%%%%%%%%%%%%%%%%%%%%%%%%%%%%%%%%%%%%%%%%%%
\section{Uniqueness of the solid angle center}
In this section, we study the uniqueness of the solid angle center of $\Omega$ of height $h$. Note that the solid angle center is not necessarily unique. Let us see in a baby case when $\dim \Omega =1$.
\begin{ex}\label{ex}{\rm Let $\Omega$ be a disjoint union of two intervals with the same length:
$\Omega = \left[-R,-1 \right] \cup \left[1,R \right]$, where $R>1$. In this case, the (solid) angle $A_\Omega^{(h)}(x)$ is given by
\[
A_{\Omega}^{(h)}(x) = \int_{\Omega}\frac{h}{(y-x)^2+h^2}dy.
\] By direct computation, we obtain the following results (see figure 2):
\begin{enumerate}
\item[(1)]If $0<h<\sqrt{R+(R+1)\sqrt{R}}$, then the (solid) angle centers of $\Omega$ are given by 
\[
x_\pm= \pm \sqrt{\sqrt{R((R+1)^2+4h^2)}-(R+h^2)}.
\]
\item[(2)]If $h\geq \sqrt{R+(R+1)\sqrt{R}}$, then $x=0$ is the unique (solid) angle center of $\Omega$.
\end{enumerate}}
\end{ex}
%%%%%%%%%%%%%%%%%%%%%%%%%%%%%%%%%%%%%%%%%%%%%%%%%%%%%%%%%%%%%%%%%%%%%%%%%%%%%%%%%%%%%%%%%%%%%%%%%%%%%%%%%%%%%%%%%
\begin{figure}[hbtp]
\centering
\scalebox{0.55}{\includegraphics[clip]{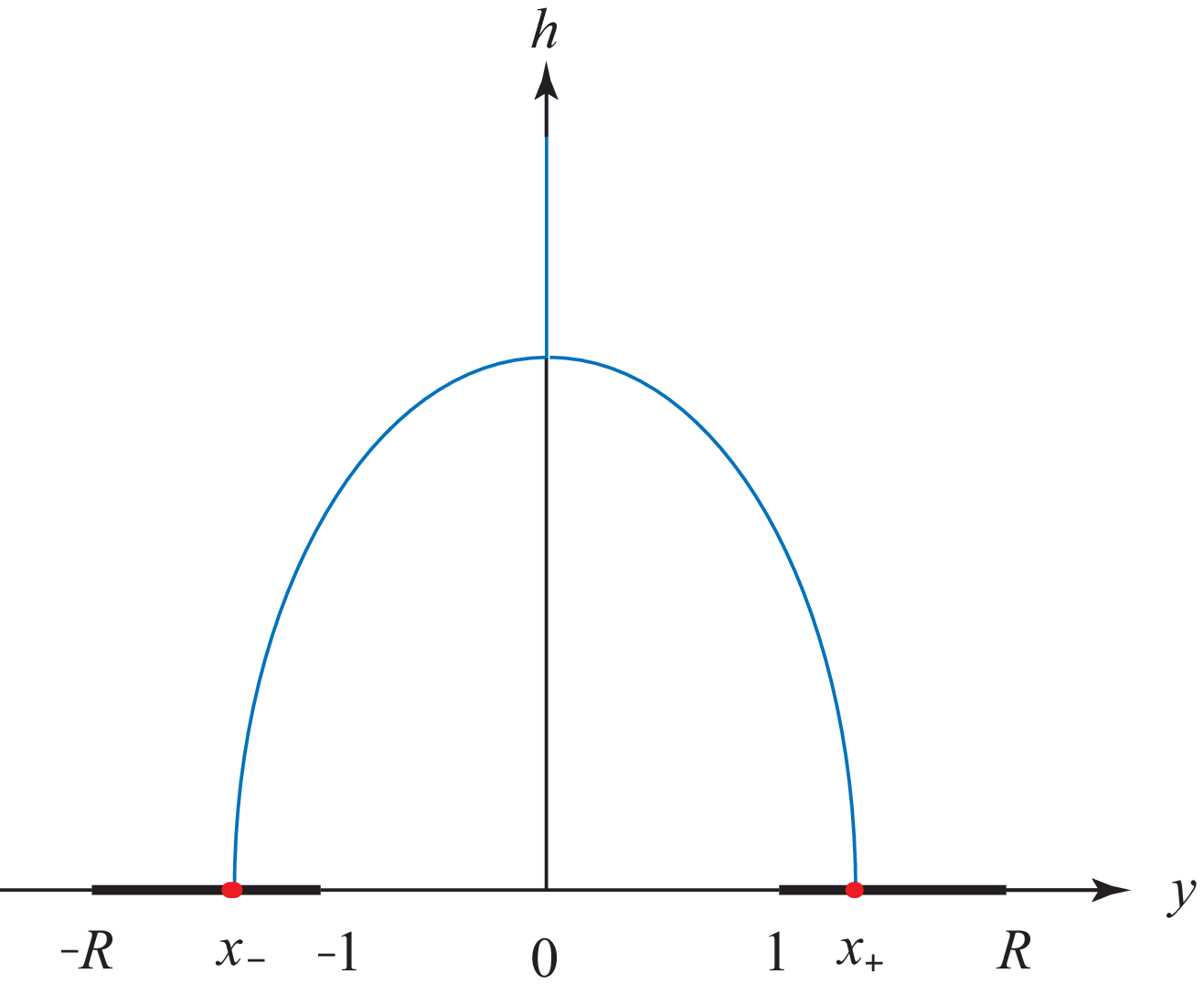}}
\caption{The locus of (solid) angle centers of height $h>0$.}
\label{example}
\end{figure}
%%%%%%%%%%%%%%%%%%%%%%%%%%%%%%%%%%%%%%%%%%%%%%%%%%%%%%%%%%%%%%%%%%%%%%%%%%%%%%%%%%%%%%%%%%%%%%%%%%%%%%%%%%%%%%%%
\subsection{Conditions for {\boldmath $h$}}
Remark that in Example \ref{ex} the solid angle center is unique if $h$ is greater than a constant determined by $\Omega$. We will show that it holds in general.

\begin{lem}\label{lem_diam} We have ${\rm diam}(\Omega) = {\rm diam}({\rm conv}(\Omega))$
\end{lem}
{\bf Proof.} There exist two points $z$ and $w$ in ${\rm conv}(\Omega)$ such that ${\rm diam}({\rm conv}(\Omega)) = \left\vert z - w \right\vert$ since ${\rm conv}(\Omega)$ is compact.\\
\hspace{1.5em} Assume that $w$ is in ${\rm conv}(\Omega)$ but not in $\Omega$. There exists two points $w'$, $w'' \in \Omega$ and a real number $0 \leq t \leq 1$ such that $w = (1-t) w' + t w''$. Then we have
\[
{\rm diam}({\rm conv}(\Omega)) = \lvert z - w \rvert < \max \left\{ \left\vert z -w' \right\vert, \left\vert z -w'' \right\vert \right\}  \leq {\rm diam}({\rm conv}(\Omega)), 
\] which is a contradiction. Hence $w \in \Omega$. With the same argument, we have $z \in \Omega$. Therefore we obtain ${\rm diam}({\rm conv}(\Omega)) \leq {\rm diam}(\Omega)$.\\
\hspace{1.5em} On the other hand, ${\rm diam}({\rm conv}(\Omega)) \geq {\rm diam}(\Omega)$ since ${\rm conv}(\Omega) \supset \Omega$. \hspace{\fill}$\Box$
\begin{lem}\label{lem_derivative<0} If $h \geq 2{\rm diam}(\Omega)$, then we have 
\[
\frac{\pd^2 A_\Omega^{(h)}}{\pd x_1^2}(x) < 0
\] for any $x \in {\rm conv}(\Omega)$.
\end{lem}
{\bf Proof.} Note that
\[
\frac{\pd ^2 A_\Omega^{(h)}}{\pd x_1 ^2}(x) = 3\int_\Omega \frac{4(y_1-x_1)^2-(y_2-x_2)^2-h^2}{(r^2+h^2)^{7/2}}dy.
\]
For any $x \in {\rm conv}(\Omega)$, if $h \geq 2{\rm diam}(\Omega)=2{\rm diam}({\rm conv}(\Omega))$, then we have
\[
\Omega \subset {\rm conv}(\Omega) \subset \left\{(y_1, y_2) \in \R^2 \vert 4(y_1-x_1)^2-(y_2-x_2)^2-h^2 \leq 0 \right\},
\] which completes the proof (see figure 3). \hspace{\fill}$\Box$\\
%%%%%%%%%%%%%%%%%%%%%%%%%%%%%%%%%%%%%%%%%%%%%%%%%%%%%%%
\begin{figure}[hbtp]
\centering
\scalebox{0.50}{\includegraphics[clip]{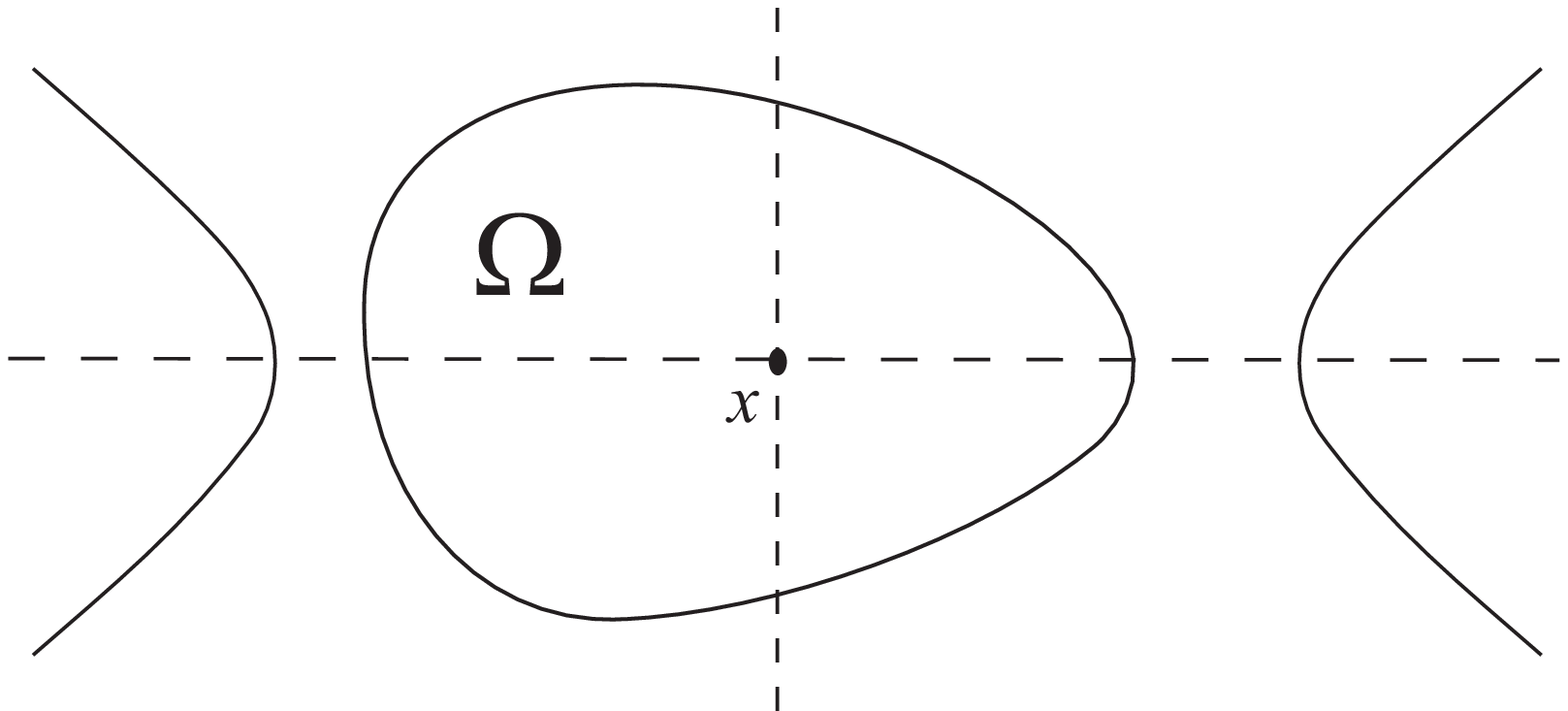}}
\caption{$\Omega \subset \left\{(y_1, y_2) \in \R^2 \vert 4(y_1-x_1)^2-(y_2-x_2)^2-h^2 \leq 0 \right\}$}
\label{condition_for_h}
\end{figure}
%%%%%%%%%%%%%%%%%%%%%%%%%%%%%%%%%%%%%%%%%%%%%%%%%%%%%%%

The radial symmetry of $\displaystyle \frac{\pd A_\Omega^{(h)}}{\pd x_1}$ implies that  the existence of two maximum points of $A_\Omega^{(h)}$ contradicts Lemma \ref{lem_derivative<0}. Hence we obtain the next corollary.
\begin{cor}\label{cor_unique} If $h \geq 2{\rm diam}(\Omega)$, then $\Omega$ has a unique solid angle center of height $h$.
\end{cor}

Let us improve the condition in Corollary \ref{cor_unique} by using the minimal unfolded region of $\Omega$. Recall that solid angle centers of $\Omega$ belong to $uf(\Omega)$.

\begin{thm}\label{thm_main1} If $h \geq 2\max\left\{\left. \left\vert z - w \right\vert \right\vert z \in uf(\Omega), w \in \pd \Omega  \right\}$, then $\Omega$ has a unique solid angle center of height $h$.
\end{thm}
{\bf Proof.} Note that we can see 
\[
\max\left\{\left. \left\vert z - w \right\vert \right\vert z \in uf(\Omega), w \in \pd \Omega  \right\}=\max\left\{\left. \left\vert z - w \right\vert \right\vert z \in uf(\Omega), w \in \pd {\rm conv}(\Omega)  \right\}
\] by the same argument as in Lemma \ref{lem_diam}. Now the theorem follows from the same argument as in Lemma \ref{lem_derivative<0} and Corollary \ref{cor_unique}. \hspace{\fill}$\Box$
%%%%%%%%%%%%%%%%%%%%%%%%%%%%%%%%%%%%%%%%%%%%%%%%%%%%%%%%%%%%%%%%%%%%%%%%%%%%%%%%%%%%%%%%%%%%%%%%%%%%%%%%%%%%%%
\subsection{Conditions for {\boldmath $\Omega$}}
We showed that the uniqueness of the solid angle center of height $h$ does not always hold for any $\Omega$ without the condition for $h$. In this subsection, we study the conditions of the uniqueness for $\Omega$. We can see the next lemma easily by Definition \ref{def_folded}.

\begin{lem}\label{lem_axial} Let $\Omega$ be given by 
\begin{equation}\label{Omega}
\Omega = \left\{ (y_1,y_2) \in \R^2 \left\vert  0 \leq y_1 \leq 1, \lvert y_2\rvert \leq f(y_1) \right. \right\}
\end{equation}
where $f:[0,1] \to \R$ is a piecewise class $C^1$ concave non-negative function. Put
\[
a=\min \left\{ \tau \in [0,1] \left\vert f(\tau)= \max_{0 \leq t \leq 1}f(t) \right. \right\} \ \
\]
and
\[
b=\max \left\{ \tau \in [0,1] \left\vert  f(\tau)= \max_{0 \leq t \leq 1}f(t) \right. \right\}.
\]
Then $uf(\Omega)$ is contained in $\displaystyle \left\{ (y_1,0) \in \Omega \left\vert  \frac{a}{2} \leq y_1 \leq \frac{1+b}{2} \right. \right\}$.
\end{lem}

Next we estimate the contribution of $\pd \Omega$ to the contour integral of the second derivative of $A_\Omega^{(h)}$.

\begin{lem}\label{lem_contribution} Let us use the notation in Lemma \ref{lem_axial} in what follows. Let $\gamma$ be an oriented curve defined by
\[
\gamma : \binom{y_1(t)}{y_2(t)} = \binom{1-t}{f(1-t)},\ \ 1-a \leq t \leq 1.
\]
Then we have 
\[
\int_{\gamma}\frac{y_1-x_1}{((y_1-x_1 )^2+y_2^2+h^2)^{5/2}}dy_2 > 0.
\] for any $\displaystyle x_1 \in \left[\frac{a}{2},{a}\right]$.
\end{lem}
{\bf Proof.}  It is equivalent to 
\[
\int_0^a\frac{(t-x_1 )f'(t)}{((t-x_1)^2+f(t)^2+h^2)^{5/2}}dt =\left(\int_0^{2x_1-a}+\int_{2x_1-a}^{x_1} + \int_{x_1}^a \right) \frac{(t-x_1 )f'(t)}{((t-x_1)^2 +f(t)^2+h^2)^{5/2}}dt <0.
\] Since $f$ is concave, we have 
\[
0 \leq f'(x_1 + \tau) \leq f'(x_1 -\tau)
\] and 
\[
0 \leq f(x_1 - \tau) \leq f(x_1 +\tau) \ \
\] for any $\tau \in \left[ 0, a-x_1 \right]$. Hence we obtain  
\begin{align*}
&\left(\int_{2x_1-a}^{x_1}+\int_{x_1}^a\right)\frac{(t-x_1)f'(t)}{((t-x_1)^2 +f(t)^2+h^2)^{5/2}}dt\\
=&\int_0^{a-x_1} \left( \frac{\tau f'(x_1 +\tau )}{(\tau ^2 +f(x_1 +\tau )^2+h^2)^{5/2}} - \frac{\tau f'(x_1 -\tau )}{(\tau ^2 +f(x_1 -\tau )^2+h^2)^{5/2}} \right) d\tau < 0
\end{align*} (see figure 4). As we have
\[
\int_0^{2x_1-a}\frac{(t-x_1)f'(t)}{((t-x_1)^2+f(t)^2+h^2)^{5/2}}dt<0,
\] it completes the proof.\hspace{\fill}$\Box$
%%%%%%%%%%%%%%%%%%%%%%%%%%%%%%%%%%%%%%%%%%%%%%%%%%%%%%%%%%%%%%%%%%%%%%%%%%%%%%%%%%%%%%%%%%%%%%%%%%%%%%%%%%%%%%%%%%%%%%
\begin{figure}[hbtp]
\centering
\scalebox{0.60}{\includegraphics[clip]{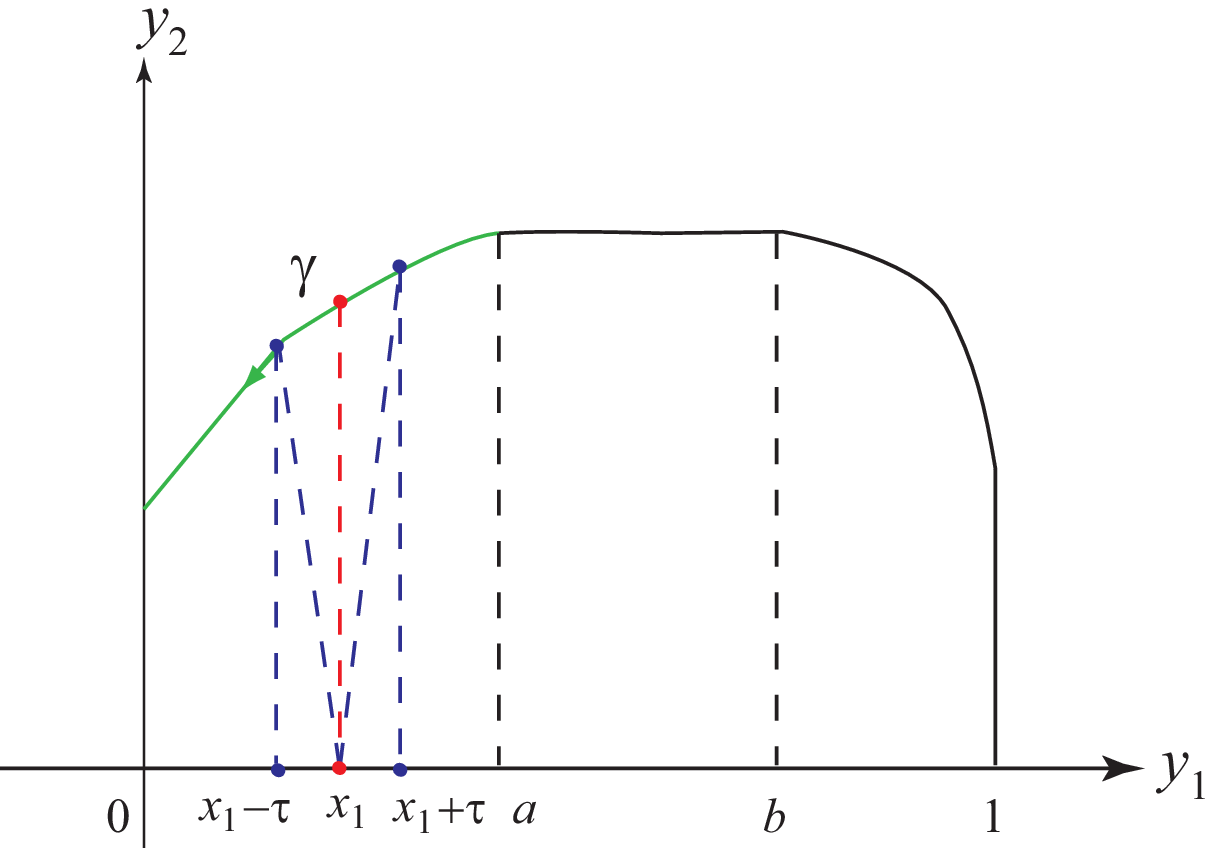}}
\caption{The path $\gamma$ and the contribution of $\pd \Omega$.}
\label{contribution}
\end{figure}
%%%%%%%%%%%%%%%%%%%%%%%%%%%%%%%%%%%%%%%%%%%%%%%%%%%%%%%%%%%%%%%%%%%%%%%%%%%%%%%%%%%%%%%%%%%%%%%%%%%%%%%%%%%%%%%%%%%%%%

We give the anti-clockwise orientation to $\pd \Omega$.

\begin{thm}\label{thm_main2} Let $\Omega$ be an axially symmetric convex body in $\R^2$ with a piecewise $C^1$ boundary. Then $\Omega$ has a unique solid angle center for any $h>0$.
\end{thm}
{\bf Proof.} By a rotation and a homothety of $\R^2$, we may assume that $\Omega$ is given by (\ref{Omega}) in Lemma \ref{lem_axial}. By Lemma \ref{lem_axial}, it suffices to show that 
\[
\frac{\pd ^2 A_{\Omega}^{(h)}}{\pd x_1^2}(x_1,0)<0
\]
for any $\displaystyle x_1 \in \left[\frac{a}{2},\frac{1+b}{2}\right]$. The formula (\ref{contour}) shows that $\displaystyle \frac{\pd ^2 A_{\Omega}^{(h)}}{\pd x_1^2}(x_1,0)$ is expressed as
\begin{equation}\label{derivative}
\frac{\pd ^2 A_{\Omega}^{(h)}}{\pd x_1^2}(x_1,0) = -3h\int_{\pd \Omega}\frac{y_1-x_1}{((y_1-x_1)^2 +y_2^2+h^2)^{5/2}}dy_2.
\end{equation} Lemma \ref{lem_contribution} implies that it is negative for $\displaystyle x_1 \in \left[ \frac{a}{2}, a \right]$. The same argument works for $\displaystyle x_1 \in \left[a, \frac{1+b}{2} \right]$. If $\displaystyle x_1 \in \left[a,b \right]$, then the right hand side of (\ref{derivative}) is obviously negative, which completes the proof. \hspace{\fill}$\Box$\\

In Theorem \ref{thm_main2} we assumed that $\Omega$ is axially symmetric. But we conjecture that the uniqueness of the solid angle center holds for any $h$ if $\Omega$ is convex. We can see a partial solution of this conjecture by using the radial function.

\begin{thm}\label{thm_main3} Suppose that $\Omega$ is convex and the minimal unfolded region of $\Omega$ is strictly contained in $\Omega$. If $0<h \leq \sqrt{2} \min \left\{ \lvert z - w \rvert \left\vert z \in uf(\Omega), w \in \pd \Omega \right. \right\}$, then $\Omega$ has a unique solid angle center of height $h$.
\end{thm}
{\bf Proof.} By using the polar coordinate, we can see 
\[
A_{\Omega}^{(h)} (x) = 2\pi - h \int_0^{2\pi} \frac{1}{\sqrt{\rho( x,\theta )^2 + h^2}} d \theta ,
\] where $\rho : \Omega \times [0,2\pi] \to \R$ is the radial function. Note that $\rho$ is well-defined on $\Omega \times [0,2\pi]$ since $\Omega$ is convex. If $0<h \leq \sqrt{2} \min \left\{ \lvert z - w \rvert \left\vert z \in uf(\Omega), w \in \pd \Omega \right. \right\}$, then for any fixed $\theta$, $\displaystyle -\frac{1}{\sqrt{\rho(x,\theta)^2 + h^2}}$ is a concave function of $x \in uf(\Omega)$. Hence the proof is completed. \hspace{\fill}$\Box$
%%%%%%%%%%%%%%%%%%%%%%%%%%%%%%%%%%%%%%%%%%%%%%%%%%%%%%%%%%%%%%%%%%%%%%%%%%%%%%%%%%%%%%%%%%%%%%%%%%%%%%%%%%%%%%%%%%%
\section{{\boldmath $r^{\al -2}$}-center}
In this section, we introduce the {\it $r^{\al -2}$-center} of $\Omega$ for $\al > 0$ from [3] and study the similar properties with solid angle function. Define the function $V_\Omega^{(\al)} : \R^2 \to \R$ for $\al > 0$ by
\[
V_\Omega^{(\al)}(x) = \left\{
\begin{array}{ll}
\displaystyle \int_\Omega r^{\al -2} dy &({\rm if} \ \  0 < \al \neq 2),\\
\displaystyle -\int_\Omega \log r dy &({\rm if} \ \  \al = 2) .
\end{array}
\right. 
\]
\begin{definition}{\rm ([3]) A point $x$ in $\R^2$ is the {\it $r^{\al -2}$-center} of $\Omega$ if it attains the maximum value of $V_\Omega^{(\al)}$ if $0<\al \leq 2$ or the minimum value of $V_\Omega^{(\al)}$ if $\al >2$.}
\end{definition}

The same arguments as in subsection 3.2 and subsection 4.2 work for the $r^{\al -2}$-center. This is because both $A_\Omega^{(h)}$ and $V_\Omega^{(\al)}$ have distance kernels which are monotonic functions of $r$. The same argument as in [4] shows the following.

\begin{thm}\label{variation2} Let $D$ be a disc in $\R^2$. If ${\rm Area}(\Omega) = {\rm Area}(D)$, then 
\begin{align*}
\max_{x \in \R^2} V_\Omega^{(\al)} (x) \leq \max_{x \in \R^2} V_D^{(\al)} (x) \ \ &({\rm if} \ \ 0< \al < 2),\\
\min_{x \in \R^2} V_D^{(\al)} (x) \leq \min_{x \in \R^2} V_\Omega^{(\al)} (x) \ \ &({\rm if} \ \ 2 \leq \al ),
\end{align*} and that equality holds if and only if $\Omega$ is a disc.
\end{thm}

Let $\Omega$ be a compact subset in $\R^2$ with non-empty interior and a piecewise $C^1$ boundary. Then it was proved in [3] (Theorem 4.6 and subsection 4.3) that there is an $r^{\al -2}$-center of $\Omega$ in the minimal unfolded region of $\Omega$ for any $\al$, and the partial derivative of $V_\Omega^{(\al)}$ can be expressed by the contour integrals as
\[
\frac{\pd V_\Omega^{(\al)}}{\pd x_1}(x) = \left\{ 
\begin{array}{ll}
\displaystyle  - \int_{\pd \Omega} r^{\al -2 } dy_2 &({\rm if} \ \ 0 < \al \neq 2),\\
\displaystyle  \int_{\pd \Omega} \log r dy_2 & ({\rm if} \ \ \al =2). 
\end{array}
\right.
\]
\begin{thm}\label{thm_main4} Let $\Omega$ be an axially symmetric convex body in $\R^2$ with a piecewise $C^1$ boundary. Then $\Omega$ has a unique $r^{\al -2}$-center for any $1<\al<3$.
\end{thm}
{\bf Proof.} By the same argument as in Lemma \ref{lem_contribution}, it holds that
\[
\int_{\gamma}\frac{y_1-x_1}{((y_1-x_1)^2+y_2^2)^{(-\al +4)/2}}dy_2 > 0
\] for any $\displaystyle x_1 \in \left[ \frac{a}{2}, a \right]$ and $1 < \al < 3$. By the same argument as the proof of Theorem \ref{thm_main2}, we obtain the conclusion. \hspace{\fill}$\Box$
\begin{rem}\label{rem_centers}{\rm Theorem \ref{thm_main3} is a partial solution of the uniqueness of $r^{\al -m }$-center for $1< \al <m+1$ which was open in [3].}
\end{rem}
%%%%%%%%%%%%%%%%%%%%%%%%%%%%%%%%%%%%%%%%%%%%%%%%%%%%%%%%%%%%%%%%%%%%%%%%%%%%%%%%%%%%%%%%%%%%%%%%%%%%%%%%%%%%%%

\no 
Department of Mathematics and Information Science,\\
Tokyo Metropolitan University,\\
1-1 Minami Osawa, Hachiouji-Shi, Tokyo 192-0397, Japan\\
E-mail: sakata-shigehiro@ed.tmu.ac.jp
\end{document}